\documentstyle{article}

\newcommand{\be}{\begin{equation}}
\newcommand{\ee}{\end{equation}}
\renewcommand{\k}[2]{\frac{#1}{#2}}

\def\s{\,\,\,\,}
\def\q{c}
\def\qq{\psi}
\def\pp{\phi}
\def\PP{\Phi}
\def\z{\zeta}
\def\a{a_0}
\def\b{\beta}
\def\ds{\dot s}
\def\kkk{(\frac{\phi}{\phi_0})^m}
\def\n{\eta}
\def\v{\lambda}
\def\ra{\rightarrow}

\title{Asymptotic Solutions of Compaction in Porous Media}

\author{Xin-She Yang      \\
  Department of Applied Mathematics and Department of Fuel and Energy \\
  University of Leeds, LEEDS LS6 9JT, UK }

\date{}

\begin{document}
\maketitle

\begin{abstract}

Compaction in reactive porous media is modelled as a  reaction-diffusion
process with a moving boundary. Asymptotic analysis is used to find
solutions for the coupled nonlinear compaction equations, and a
traveling wave solution is obtained above the reaction zone. \\


\noindent {\it Citation detail:} X. S. Yang, Asymptotic solutions of compaction in porous media,
{\it Applied Mathematics Letters}, {\bf 14}, 765-768 (2001).

\end{abstract}

\section{Introduction}

The accurate modelling of compaction and reactive flow in porous media
such as sand and shale sediments is very important in civil engineering
and oil industry [1]. The general mathematical model of compaction and
reactive mineral flow considers the fluid-sediment system as a porous
medium consisting of multiple mineral species. The interstitial volume
of the porous solid phase is saturated with pore fluids. Due to the action
of gravitational overburden loading and the density difference between
the two phases, the solid phase compacts by reducing the porosity
(volume fraction of the pores), thus leading to
the expulsion of the pore fluid out of porous media.
During compaction and continuous burial, the mineral species
such as water-rich clay smectite  react with pore fluids and are then
transformed into a more stable mineral species such as illite, releasing
free water into the porous environment[1,2]. In this paper,
a reaction-diffusion model together with some asymptotic
analysis is presented.

\section{Mathematical  Model}

Let the volume fractions of solid reactant species (smectite) and
water be $\qq,\, \pp$, respectively. By proper
non-dimensionalization and  appropriate scalings [3,4,5] in
a 1-D basin $0<z<h(t)$, where $h(t)$ is the ocean floor and
$z=0$ is the basement rock, we can
write down the non-dimensional compaction equations as
\be
\qq_t=-e^{\b (h-z-z^*)} \qq -\k{ \v}{1-\pp_0} [\qq \kkk (\pp_z-\pp)]_z.
\label{equ-100}
\ee
\be
\pp_t=\v[ \kkk (\pp_{z}-\pp)]_z + \k{\a}{\b} e^{\b (h-z-z^*)} \qq,
\s  m \ge 7.  \label{equ-200}
\ee
The boundary conditions are
\be
\pp_z-\pp=0, \s {\rm at } \s z=0,
\label{bbb-1}
\ee
\[ \pp=\phi_{0}, \s, \qq=\qq_0, \]
\be
\dot h(t)=\dot s +\k{\v}{1-\pp_0} \kkk (\pp_z-\pp) \,\,
{\rm at} \,\,z=h(t),
\label{bbb-2}
\ee
where $\v=O(1)$ and $\b \gg 1$ are compaction constant and
non-dimensional activation energy of the one-step dehydration
mineral reaction with a critical temperature $z^*$ and the
amount $\a=O(1)$ of free water released from the reaction.
$\dot s$ is the rate of new sediment accumulation at the
basin top, and thus can be taken as a prescribed constant
($\dot s=1$). $\pp_0$ and $\qq_0$ are the initial values of
$\pp$ and $\qq$ at the ocean floor $z=h(t)$, respectively.
In addition, all the variables ($\pp, \, \qq, \, z, \, t$) and
the parameters ($\v, \, m, \, \pp_0, \qq_0, \,\dot s$) are non-negative.
The reaction term $\exp[\b (h-z-z^*)]$ is only dominant within
a thin region of a width of $O(1/\b)$ near $h-z-z^* \approx 0$,
in other words, the reaction region is located at $z \approx \theta^*$
defined as
\be
\theta^*=h-z^*.
\ee
Clearly, the present problem is a non-linear diffusion problem
with a  boundary $h(t)$ moving at a speed of $\dot h(t)$,
which can be solved numerically by using the predictor/corrector
implicit finite-difference method.

\section{Travelling Wave Solution}

Numerical simulations [3,6]
and real data observations [7] imply that the moving
boundary  $h(t)$ moving essentially at a nearly constant speed
$\dot h=\q$, although the specific value $\q$ depends on the boundary
conditions, and is yet to be determined.
The constant moving boundary implies there exist travelling wave
solutions. Thus, we define a new variable by
\be
\z=z-\theta^*=z-h(t)+z^*, \s -\theta^* \le \z \le z^*,
\ee
so that the model equations (\ref{equ-100}) and
(\ref{equ-200}) become
\be
-\q \qq_{\z}=-e^{-\b \z} \qq -\k{\v}{1-\pp_0} [\qq \kkk (\pp_{\z}-\pp)]_{\z}.
\label{equ-500-1}
\ee
\be
-\q \pp_{\z}=\v [\kkk (\pp_{\z}-\pp)]_{\z} + \k{\a}{\b} e^{-\b \z} \qq,
\label{equ-500-2}
\ee
The fact that $\b \sim m \gg 1$ allows us to seek asymptotic
solutions in different regions with these distinguished limits. $\b \gg
1$ implies that the reaction is essentially restraint in a very
narrow zone with a width of $O(1/\b)$ at $z=\theta^*$. Above this region
($z > \theta^*$), we have $\z > 0$ so that $\exp(-\b \z) \ra 0$.
Below this reaction region ($\z < 0$), the reaction is essentially
completed (i.e., the volume fraction of smectite $\qq \ra 0$)
and consequently $\qq \exp(-\b \z) \ra 0$. We shall
see that this is true below in equation (\ref{psi-solution}) because
$\qq \ra 0$ as $\z \ra -\infty $ (or $\n \ra -\infty$).

\subsection*{Outer Solutions}

In the outer region above the reaction zone ($\z >0 $),
the reaction terms are negligible, then we have
\be
\q \qq_{\z}=\k{\v}{1-\pp_0} [\qq \kkk (\pp_{\z}-\pp)]_{\z}, \label{equ-300}
\ee
\be
-\q \pp_{\z}=\v [\kkk (\pp_{\z}-\pp)]_{\z}. \label{equ-400}
\ee
The integrations together with top boundary conditions (\ref{bbb-2})
give
\be
\qq =\k{\ds \qq_0}{\q  -\k{\v}{1-\pp_0} \kkk (\pp_{\z}-\pp)},
\ee
and
\be
\q \pp+\v \kkk (\pp_{\z}-\pp)=\q \pp_0+(\q -\ds)(1-\pp_0), \label{solution-1}
\ee
whose solution can be further written as a quadrature. This solution is
only  valid in the region above the reaction zone ($z=\theta^*$).
On the other hand,
the travelling solution will not be appropriate in the region below the
reaction zone ($\z < 0$) because the exponential term
$(\pp/\pp_0)^m \ll 1$ due to $\pp < \pp_0$ and $m \gg 1$. However,
we can define a typical value of $\pp^*$ by
\be
\pp^*=\pp_0 e^{-\k{\ln m }{m}},
\ee
so that $\pp \sim \pp^*$ in the region below the reaction zone ($\z < 0$).
The reaction term is still negligible. Thus we rewrite
equations (\ref{equ-100}) and (\ref{equ-200}) in terms of $\PP$ defined by
\be
\pp=\pp^* e^{\k{\PP}{m}}=\pp_0 e^{\k{\PP - \ln m}{m}},
\ee
then we have
\be
\qq_{t}=-\k{\pp^* \v}{m (1-\pp_0)}[\qq e^{\PP}
(\k{1}{m} \PP_{z}-1)]_{z}.    \label{equ-1-1}
\ee
\be
\pp^* \PP_{t}=\v \pp^* [e^{\PP} (\k{1}{m} \PP_{z}-1)]_{z},
\label{equ-2-2}
\ee
By using $1/m \ll 1$, the above equations becomes
\be
\qq_{t} \approx 0,
\ee
\be
\PP_{t}+\v e^{\PP} \PP_{z}=0.
\ee
It is straightforward to write down the solution together with the boundary
condition $\PP_z=m$ at $z=0$. We have
\be
\PP =\ln ( \k{1+m z}{1+m \v t} ). \label{equ-600}
\ee
As $z=h(t)$, we have $\PP_{\infty} =\ln \{(1+m h)/(1+m \v t)\}
\ra \ln (h/\v t)=\ln (\q/\v)$ as $t \ra \infty$ or $ m h(t) \ra \infty$
due to $d h/d t=c=const$ and $h=0$ at $t=0$. However, $z=h(t)$ is not
usually reached since solution (\ref{equ-600}) is below the reaction region.
When $z \gg 1$, the above solution shall match the inner solutions as $\n
\ra -\infty$.

\subsection*{Inner Solutions}

In the reaction zone, we use the stretched variables defined by
\be
\n=\b \z+\ln \b,  \s \pp=\pp^* e^{\k{\PP}{m}}=\pp_0 e^{\k{\PP - \ln m}{m}},
\ee
so that equations (\ref{equ-500-1}) and (\ref{equ-500-2}) become
\be
-\q \qq_{\n}=-e^{-\n} \qq -\k{A \v \pp^*}{\b (1-\pp_0)}[\qq e^{\PP}
(A \PP_{\n}-1)]_{\n}.    \label{equ-500}
\ee
\be
-\q \pp^* \PP_{\n}=\v \pp^* [e^{\PP} ( A \PP_{\n}-1)]_{\n} + \k{\a}{A}
e^{-\n} \qq, \label{equ-600}
\ee
where $A=\k{\b}{m}=O(1)$. By using $1/\b \ll 1$,  equation
(\ref{equ-500}) becomes
\be
\q \qq_{\n}=e^{-\n} \qq,
\ee
whose solution is
\be
\qq =C \exp[-\k{1}{\q}
e^{-\n}], \s C=\qq_0 \exp[\k{1}{c} e^{(\b z^*-\ln \b)}],
\label{psi-solution}
\ee
where $C$ depends on $\q$, and $\q$ will be determined later in
(\ref{match-equ}) by matching.
Substituting this solution into (\ref{equ-600}), integrating once
from $-\infty$ to $\n$ and using $\PP \ra \PP_{\infty}$ as $\n \ra
-\infty$,  we get
\be
\q \pp^* \PP + \v \pp^* e^{\PP} (A \PP_{\n}-1)-B= -\k{\q \a}{A} C,
\label{solution-2}
\ee
where $B=\q \pp^* \PP_{\infty}-\v \pp^* \exp(\PP_\infty)$. As $\n \ra
\infty$, we have
\be
[\q \pp^* \PP + \v \pp^* e^{\PP}(A \PP_{\n}-1)]^{+\infty}_{-\infty}=
- \k{\q \a}{A} C,
\ee
which implies a jump through the reaction region.

\subsection*{Matching}

By rewriting solution (\ref{solution-2}) in terms of $\pp$,
we have approximately
\be
\q \pp +\v \kkk (\pp_{\z}-\pp) \approx
\q \pp^* \PP_{\infty}-\v \pp^* e^{\PP_{\infty}} -C \k{\q \a}{A},
\ee
and matching this to solution (\ref{solution-1}), we have
\be
\q \pp_0 +(\q-\dot s)(1-\pp_0)=\q \pp^* \PP_{\infty}-\v \pp^* e^{\PP_{\infty}}
-C \k{\q \a}{A},
\ee
which determines $\q$, leading to
\be
\q=\k{\ds (1-\pp_0)}{1-\pp^* \PP_{\infty} + \k{\a C}{A}} -
\k{\v \pp^* e^{\PP_{\infty}}}{(1-\pp^* \PP_{\infty} + \k{\a C}{A})},
\label{match-equ}
\ee
Since $\PP_{\infty}$ is a function of $\q$, we now have an implicit equation
for $\q$, which depends essentially on the initial values of $\pp_0$ and
$\qq_0$.

In summary, although it is very difficult to seek directly  solutions
for the couple nonlinear reaction-diffusion equations, we use a hybrid method
to get the asymptotic solutions and travelling wave solutions in
different regions. The matching of these solutions can thus determine
the boundary moving velocity $\dot h(t)=c$, which shows how the reaction
inside the porous media affect the evolution of its top boundary.

\section*{References}

\begin{description}

\item[1]   D. Eberl and J. Hower, Kinetics of illite formation, {\it Geol.
Soc. Am. Bull.},  {\bf 87}(10) 1326-1330(1976).

\item[2]   H. J. Abercrombie,  {\it et al.}, Silica activity and the smectite-illite reaction, {\it Geology}, {\bf 22}, 539-542 (1994).

\item[3]   X. S. Yang, {\it Mathematical modelling of compaction
           and diagenesis in sedimentary basins},
           PhD dissertation, Oxford University (1997).

\item[4]   D. M. Audet and A. C. Fowler, A mathematical model for
compaction in sedimentary basins, {\it Geophys. Jour. Int.}, {\bf 110}(3),
577-590 (1992).

\item[5]   A. C. Fowler and X. S. Yang, Fast and slow compaction
in sedimentary basins, {\it SIAM J. Appl. Math.}, 1998, {\bf 59}(1),
365-385 (1998).

\item[6]   X. S. Yang, Nonlinear viscoelastic compaction in sedimentary
           basins, {\it  Nonlinear Proc. Geophys.},{\bf 7}, 1-8 (2000).

\item[7]   I. Lerche, I.  Basin analysis: quantitative methods, Vol. I,
Academic Press, San Diego, California (1990).

\end{description}

\end{document}